\documentclass[11pt]{amsart}
\usepackage{amsmath}

\input xy
\xyoption{all}
\newcommand{\Proj}{{\mathbb P}}
\newcommand{\R}{{\mathbb R}}

\newcommand{\vol}{\mbox{\it vol}}
\renewcommand{\S}{{\mathcal S}}

\newcommand{\leg}{{\mathcal L}}
\newcommand{\Proof}{\noindent{\it Proof.\ }}

\newtheorem{theorem}{Theorem}[section]
\newtheorem{corollary}{Corollary}[section]
\newtheorem{lemma}{Lemma}[section]
\newtheorem{proposition}{Proposition}[section]
\newtheorem*{conjecture}{Conjecture}
\newtheorem*{theorem*}{Theorem}

\theoremstyle{definition}
\newtheorem{definition}{Definition}[section]

\title{Dual Spheres Have the Same Girth}

\author{J.C. \'Alvarez Paiva}

\thanks{This work was partially funded by a {\it cr\'edit aux
chercheurs\/} from the FNRS, Belgium}

\address{J.C. \'Alvarez Paiva, Polytechnic University, 6 MetroTech
Center, Brooklyn, New York, 11201.}

\email{jalvarez@duke.poly.edu}



\begin{document}%

\begin{abstract}
Symplectic and Finsler geometry are used to settle a conjecture
of Sch\"affer stating that the {\it girth\/} of a normed space --- the infimum
of the lengths of all closed, rectifiable, centrally symmetric curves on its
unit sphere --- equals the girth of its dual.
\end{abstract}

\maketitle

\section{Introduction}

A remark that deserves to be better known is that the projectivized of
a normed space inherits a natural metric that plays an important role
in the study of its geometric properties. If the distance between two
points on the unit sphere $\S \subset X$ is defined as the infimum of
the length of all rectifiable curves on $\S$ that join them, this distance,
being invariant under the antipodal map, induces a metric on the projectivized
space $\Proj X$. Since, in general, this metric is a low-regularity
Finsler metric, this functorial construction opens the way for the
use of metric and differential geometry in the study of affine
invariants of convex bodies and the geometry of normed spaces.
Indeed, any metric or Finsler invariant of $\Proj X$ is an affine
invariant of the convex hypersurface $\S$ and an isometry
invariant of $X$.

An specially interesting invariant of $\Proj X$ is its
{\it systole:\/} the infimum of the lengths of all non-contractible
rectifiable curves. The resulting invariant of isometry classes of
normed spaces was extensively studied in \cite{Schaffer:2} by
J.~J.~Sch\"affer, who used the term {\it girth\/} for twice the
systole of $\Proj X$. Among other results, Sch\"affer found estimates
for the girth of an $n$-dimensional normed space, showed that the
girth was a continuous (and even Lipschitz) function in the
Banach-Mazur compactum of isometry classes of $n$-dimensional normed
spaces, and posed two beautiful conjectures:

\begin{conjecture}
The girth of a normed space equals the girth of its dual.
\end{conjecture}

\begin{conjecture}
The girth of a normed space of dimension greater than two is less
than or equal to $2\pi$. Equality holds if and only if the norm comes
from an inner product.
\end{conjecture}

This paper settles the first of these conjectures by presenting a
``proof from the book" that uses only some elementary symplectic and
Finsler geometry.

When $X$ is finite-dimensional and the metric induced on $\S$, and
$\Proj X$, is a smooth Finsler metric, we shall also see that the length
spectrum  --- the set of lengths of all periodic geodesics --- of
$\S$ equals the length spectrum of $\S^* \subset X^*$ and that the
length spectrum of $\Proj X$ equals the length spectrum of $\Proj X^*$.
The paper ends by showing that the Cauchy-Crofton formula on normed
spaces also follows from the same symplectic  constructions.

\section{Short review of symplectic and Finsler geometry}

Throughout the paper we shall be studying the geometry of a certain
class of hypersurfaces in cotangent bundles. While many of the definitions
and notions reviewed below are familiar to symplectic geometers,
I have included a short exposition as a courtesy to readers with
other backgrounds and in order to fix notation. The references for
this section are the books \cite{McDuff-Salamon} and
\cite{Arnold-Givental}, as well as the thesis of Cieliebak,
\cite{Cieliebak}.

Let $M$ be an $n$-dimensional smooth manifold and let
$\pi : T^*M \rightarrow M$ be its cotangent bundle. It is well known that
$T^*M$ carries a tautological or {\it canonical\/} 1-form, $\alpha$, whose
value at a vector $\dot{\xi} \in T_{\xi} T^*M$ is $\xi(D\pi(\dot{\xi}))$. The
two-form $\omega := - d\alpha$ is the canonical or {\it symplectic\/} two-form.

A compact hypersurface $H \subset T^*M$ is said to be of {\it contact-type\/}
if the pullback of $\alpha$ to $H$ is a contact form (i.e., the pullback of
$\alpha \wedge (d\alpha)^{n-1}$ to $H$ never vanishes).
The hypersurfaces of contact-type that we shall consider in this paper are
all of the following form: $M$ is a compact manifold and for every $x \in M$
the intersection $H \cap T_{x}^* M$ is a convex hypersurface of $T_{x}^* M$
enclosing the origin.

By a slight abuse of notation, if $H \subset T^*M$ is a hypersurface of
contact-type, the pullback of $\alpha$ to $H$ will still be denoted as
$\alpha$. The geometry of $(H,\alpha)$ is quite rich:
first of all $H$ can be given an orientation by specifying that
$\alpha \wedge (d\alpha)^{n-1}$ be a (positive) volume form. We can then
define the {\it volume\/} of $(H,\alpha)$ by the formula
$$
\vol(H,\alpha) := \frac{1}{n!} \int_{H} \alpha \wedge (d\alpha)^{n-1} .
$$
The integral of $\alpha$ along a curve $\sigma :[a,b] \rightarrow H$ is
called the {\it action\/} of $\sigma$.

The fact that the form $\alpha \wedge (d\alpha)^{n-1}$ is a volume form is
equivalent to the following statement:
for each $\xi \in H$ the kernel of the form
$\omega_{\xi} = -d\alpha_{|T_{\xi}H}$,
$$
\hbox{Ker } \omega_{\xi} := \{\dot{\xi} \in T_{\xi} H :
                            \omega_{\xi}(\dot{\xi},\cdot) = 0 \} \, ,
$$
is a one-dimensional subspace. Moreover, this line can be oriented by
saying that $\dot{\xi}$ points in the positive direction if
$\alpha(\dot{\xi}) > 0$. The line field $\xi \mapsto \hbox{Ker } \omega_{\xi}$
is called the {\it characteristic line field.} Its positively oriented integral
curves are the {\it characteristics\/} of $(H,\alpha)$.

The set of the actions of all periodic characteristics is called the
{\it action spectrum.} Note that since
a periodic characteristic can be traversed $k$ times, $k$ a positive integer,
if a number $a$ is in the action spectrum, then so is $k a$.

\begin{definition}
Let $(H_{1},\alpha_{1})$ and $(H_{2},\alpha_{2})$ be two hypersurfaces of
contact-type. A diffeomorphism $\phi : H_{1} \rightarrow H_{2}$ is said to
be an {\it exact contactomorphism\/} if either $\phi^* \alpha_{2} - \alpha_{1}$
or $\phi^* \alpha_{2} + \alpha_{1}$ is an exact 1-form.
\end{definition}

\begin{proposition}\label{invariants}
Two exact contactomorphic hypersurfaces of contact-type have the same volume and
action spectrum.
\end{proposition}

\Proof
Let $\phi$ be a diffeomorphism between two $(2n-1)$-dimensional hypersurfaces
of contact-type $(H_1,\alpha_1)$ and $(H_2, \alpha_2)$ that satisfies
$\phi^* \alpha_2 = -\alpha_1 + df$, where $f$ is a smooth function on $H_1$.

The equality $\vol(H_1,\alpha_1) = \vol(H_2,\alpha_2)$ follows from the fact that
$H_1$ is closed and that $\phi^* \, \alpha_2 \wedge (d\alpha_2)^{n-1}$ equals
$(-1)^{n} \, \alpha_2 \wedge (d\alpha_2)^{n-1}$ plus an exact form.

It is also easy to see that $\phi$ maps the characteristic line field in
$(H_1,\alpha_1)$ to the characteristic line field in $(H_2,\alpha_2)$. However,
since we defined characteristics to be positively oriented, we cannot say that
$\phi$ maps characteristics to characteristics. What we can say is that
if $\sigma : [a,b] \rightarrow H_1$ is a closed characteristic, then the curve
$t \mapsto \phi(\sigma(-t))$ is a closed characteristic with the same action.

This is enough to show that the action spectra of $(H_1,\alpha_1)$ and
$(H_2,\alpha_2)$ coincide.

If the map $\phi$ satisfies $\phi^* \alpha_2 = \alpha_1 + df$, then it does
take closed characteristics to closed characteristics with the same action, and
the proof is only easier.
\qed

For the rest of the paper all hypersurfaces of contact-type that we shall consider
will be unit co-sphere bundles of Finsler metrics, reversible and non-reversible.
Before reviewing Finsler metrics, let us describe these hypersurfaces in simple
terms:

\begin{definition}
A smooth hypersurface in a finite-dimensional real vector space $V$ is said
to be {\it quadratically convex\/} if its principal curvatures are positive
for any auxiliary Euclidean structure on $V$. Equivalently, a smooth hypersurface
is quadratically convex if and only if its osculating quadrics are
all ellipsoids.
\end{definition}

A hypersurface $H \subset T^* M$ is said to be {\it optical\/} if in each
cotangent space $T_{x}^* M$ the intersection $H \cap T_{x}^* M$ is a
quadratically convex hypersurface enclosing the origin. When the manifold
$M$ is compact, optical hypersurfaces of $T^* M$ are of contact-type.

To make the link between optical hypersurfaces and Finsler metrics, we extend
slightly the notion of duality in normed spaces:

Let $V$ be a finite-dimensional real vector space and let $\S \subset V$ be a
quadratically convex hypersurface enclosing the origin. If $q$ is a point in $\S$,
there is a unique covector $\xi \in V^*$ such that the hyperplane
$\xi = 1$ is tangent to $\S$ at $q$ and the half-space $\xi \leq 1$ contains
$\S$. The map that sends $q$ to $\xi$ is called the {\it Legendre transform\/}
and will be denoted by
$$
\leg : \S \longrightarrow V^* .
$$
The image of $S$ under $\leg$, the {\it dual\/} of $\S$, is again a quadratically
convex hypersurface that encloses the origin. It will be denoted by $\S^*$. It
is well-known, and easy to verify, that $\S$ is also the dual of $\S^*$.

\begin{definition}
Let $M$ be a smooth manifold and let $TM \setminus 0$ denote its
tangent bundle with the zero section deleted. A {\it Finsler
metric\/} on $M$ is a smooth function
$$
\varphi : TM \setminus 0 \longrightarrow \R
$$
that is positively homogeneous of degree one and such that on every tangent
space $T_x M$ the hypersurface
$$
S_{x}M := \{v \in T_x M : \varphi(v) = 1 \}
$$
is quadratically convex and encloses the origin. The Finsler metric will
be called {\it reversible\/} if $S_x M$ is symmetric with respect to the
origin (i.e., if and only if $\varphi(v) = \varphi(-v)$).
\end{definition}

We follow the standard terminology in Finsler geometry and define a
{\it Minkowski space\/} to be a normed space whose unit sphere is
quadratically convex. Clearly, any submanifold of a Minkowski space inherits
a natural reversible Finsler metric. In this paper we are mainly concerned
with the geometry of the Finsler metric inherited by the unit sphere of a
Minkowski space.

A reversible Finsler metric $\varphi$ on a manifold $M$ allows us to define
the length of a smooth curve $\gamma : [a,b] \rightarrow M$ by the
equation
$$
\mbox{length of } \gamma := \int_{a}^{b} \varphi(\dot{\gamma}(t)) \, dt .
$$
Using this, we can define the distance between two points $x$ and
$y$ in $M$ as the infimum of the lengths of all smooth curves
joining $x$ and $y$. As usual, a geodesic on a Finsler manifold is a curve
that locally minimizes length. The quadratic convexity of the tangent unit
spheres guarantees the local uniqueness of geodesic segments. When the
Finsler metric is non-reversible, we must settle for oriented versions
of the notions of length, distance, and geodesic.

Notice that the restriction of Sch\"affer's conjecture for Minkowski spaces
states that {\it the length of the shortest centrally symmetric closed geodesic
on the unit sphere of a Minkowski space equals the length of the shortest
centrally symmetric closed geodesic on the unit sphere of its dual.}

A Finsler metric on a manifold $M$ defines an optical hypersurface on its
cotangent bundle: if  $S^*_{x}M \subset T_{x}^* M$ denotes the dual of $S_x M$,
then the {\it unit co-sphere bundle\/} of the Finsler manifold $(M,\varphi)$,
$$
S^* M := \bigcup_{x \in M} S^{*}_{x} M \subset T^* M \, ,
$$
is an optical hypersurface. Many of the geometric quantities and objects
on $(M,\varphi)$ can be read from the geometry of $(S^* M, \alpha)$.

Let us denote by $\leg : SM \rightarrow S^*M$ the map that takes a unit vector
$v \in S_{x}M$ and sends it to $\leg_{x}(v)$, where
$\leg_x : S_{x}M \rightarrow T_{x}^* M$ is the Legendre transform at the point $x$.
If $\gamma$ is a smooth curve on $M$ parameterized by arc-length, the curve
$\leg \circ \dot{\gamma}$ is a smooth curve on $S^* M$, and it is easy to see that
the length of $\gamma$ equals the action of $\leg \circ \dot{\gamma}$:
\begin{equation} \label{action}
\mbox{ length of } \gamma := \int_{\gamma} \varphi =
                             \int_{\leg \circ \, \dot{\gamma}} \alpha .
\end{equation}

It is also well-known that {\it a curve $\gamma : [a,b] \rightarrow M$ is a
geodesic if and only if the curve $\leg \circ \dot{\gamma}$ is a
characteristic of $(S^* M,\alpha)$.} This, together with equation (1), shows
that the {\it length spectrum\/} of a Finsler manifold --- the set of lengths
of its periodic geodesics --- equals the action spectrum of its unit co-sphere
bundle.

We may also define the volume of an $n$-dimensional Finsler manifold $M$ as the
volume of $(S^* M, \alpha)$ divided by the volume of the
Euclidean $n$-dimensional unit ball. This definition was first introduced from a
convex-geometric viewpoint by Holmes and Thompson (see \cite{Holmes-Thompson}
and \cite{Thompson}), and is usually referred to as the
{\it Holmes-Thompson volume\/} of $M$.

If we define the {\em unit co-disc bundle\/}, $D^*M$, of a Finsler manifold
$M$ to be the open set bounded by $S^*M$ and containing the zero section of
$T^*M$, then, by Stokes theorem, the Holmes-Thompson volume of $M$
equals the symplectic volume of $D^* M$ divided by the volume of the Euclidean
$n$-dimensional unit ball.

\section{A geometric construction}

The solution of Sch\"affer's conjecture and the other results in this paper
depend on a simple geometric construction:

Let $V$ be a finite-dimensional real vector space and let $\S_{1} \subset V$
and $S_{2}^* \subset V^*$ be two quadratically convex hypersurfaces.  If
$q \in \S_{1}$ and $P \in V^*$, let us view $T_{q} \S_{1}$ as a vector subspace
of $V$ and define $P_q \in T_{q}^* \S_1$ to be the restriction of $P$ to
$T_{q} \S_{1}$. Likewise, if $P \in \S_{2}^*$ and $q \in V$, the restriction of
$q$ to $T_{P}\S_2^*$ will be denoted by $q_{P}$.

Let us now consider the sets
\begin{eqnarray*}
\overline{D}^* (\S_{1},\S_{2}^*) & := &
 \{P_q \in T_{q}^* \S_{1} : q \in \S_1 \, , \, P \in \S_{2}^* \}
  \subset T^* \S_1 , \\
\overline{D}^* (\S_{2}^*,\S_{1}) & := &
 \{q_P \in T_{P}^* \S_{2}^* : q \in \S_1 \, , \, P \in \S_{2}^* \}
  \subset T^* \S_{2}^*
\end{eqnarray*}
together with their interiors, $D^*(\S_{1},\S_{2}^*)$ and
$D^*(\S_{2}^*,\S_{1})$, and their boundaries, $S^*(\S_{1},\S_{2}^*)$
and $S^*(\S_{2}^*,\S_{1})$.

\begin{lemma}\label{contact-lemma}
The hypersurfaces $S^*(\S_{1},\S_{2}^*) \subset T^* \S_1$ and
$S^*(\S_{2}^*,\S_{1}) \subset T^* \S_{2}^*$ are exact contactomorphic.
\end{lemma}

\Proof
The proof shall consist of two parts: constructing a diffeomorphism
$\phi$ between $S^*(\S_{1},\S_{2}^*)$ and $S^*(\S_{2}^*,\S_{1})$, and
showing that $\phi$ is an exact contactomorphism.

\medskip \noindent
{\it Construction of $\phi$:} First we remark that if
$(q,p) \in S^*(\S_{1},\S_{2}^*)$, then there exists a unique
$P \in \S_{2}^*$ such that $P_q = p$. Indeed, the quadratic convexity of
$\S_{2}^*$ implies that the projection $P \mapsto P_q$ defines a diffeomorphism
between its singular set, {\it the shadow boundary\/} of the projection,
and the boundary of its image in $T_q^* \S_1$.

If $(q,p) \in S^*(\S_{1},\S_{2}^*)$ and $P$ is the unique point in $\S_2^*$
such that $P_q = p$, we define $\phi(q,p) := (P,q_P) \in T_{P} \S_{2}^*$.
Note that it is not yet obvious that $\phi(q,p)$ belongs to
$S^*(\S_{2}^*,\S_{1})$.

That this is true follows at once from the following general statement:
{\it a point $P \in \S_2^*$ is in the shadow boundary of the projection
$X \mapsto X_q$, $q \in \S_{1}$, from $\S_2^*$ to $T_{q}^* \S_1$ if and
only if the point $q$ is in the shadow boundary of the projection
$x \mapsto x_P$ from $\S_1$ to $T_P^* \S_2^*$.}

Indeed, the condition that $P$ be in the shadow boundary of the projection
$X \mapsto X_q$ is equivalent to the degeneracy of the bilinear form given by
the restriction of the dual pairing in $V^* \times V$ to
$T_P \S_2^* \times T_q \S_1$. This is the same condition for $q$ to be
in the shadow boundary of the projection $x \mapsto x_P$.

The map $\phi$ is obviously invertible, and the quadratic convexity of both
$\S_1$ and $\S_2^*$ implies that both it and its inverse are smooth.

\medskip \noindent
{\it The map $\phi$ is an exact contactomorphism:}
To see that  $\phi^* \alpha_2 + \alpha_1$ is an exact $1$-form, we borrow an
old trick from classical mechanics and consider the graph of $\phi$ as a
submanifold of $S^*(\S_{1},\S_{2}^*) \times S^*(\S_{2}^*,\S_{1})$:
$$
\Gamma_\phi := \{(q,p;P,Q) \in S^*(\S_{1},\S_{2}^*) \times S^*(\S_{2}^*,\S_{1}) :
             q_P = Q \mbox{ and } P_q = p \} \, .
$$
We may consider the $1$-form
$$
\alpha_{2} + \alpha_{1} = Q \cdot dP + p \cdot dq
$$
as a form on $S^*(\S_{1},\S_{2}^*) \times S^*(\S_{2}^*,\S_{1})$. The pullback of
this form to $\Gamma_\phi$ becomes
$$
q \cdot dP + P \cdot dq = d(P \cdot q) ,
$$
and, therefore, $\phi$ is an exact contactomorphism.
\qed

\begin{lemma}\label{symplectic-lemma}
There exists an anti-symplectic diffeomorphism between the open sets
$D^*(\S_{1},\S_{2}^*) \subset T^* \S_1$ and
$D^*(\S_{2}^*,\S_{1}) \subset T^* \S_{2}^* $.
\end{lemma}

\Proof
The idea of the proof is basically the same as that of the previous lemma.
There is only the slight complication that if $p \in T_q^* \S_1$ is in
$D^*(\S_1,\S_2^*)$, then there exist two points $P^+$ and $P^-$ in $\S_2^*$
such that $P_q^+ = P_q^- = p$. To distinguish between these points and define
a diffeomorphism $\Phi$ between $D^*(\S_1,\S_2^*)$ and $D^*(\S_2^*,\S_1)$, we
need to play with the orientations.

Let us fix an orientation on $V$. This induces an orientation on $V^*$ and,
using the natural co-orientation of the hypersurfaces $\S_1$ and $\S_2^*$,
on $\S_1$ and $\S_2^*$. In turn, the orientations on $\S_1$ and $\S_2^*$
induce orientations on each of their cotangent spaces.

If $(q,p) \in D^*(\S_1,\S_2^*)$, there is a unique point $P \in \S_2^*$
such that $P_q = p$ and such that the differential of the projection
$X \mapsto X_q$ from $\S_2^*$ to $T_q^* \S_1$ at the point $P$ is
orientation preserving. We define $\Phi(q,p) := (P,q_P)$.

Exactly as in the proof of the previous lemma, one shows that
$\Phi^*\alpha_2 + \alpha_1$ is an exact $1$-form, and
$\Phi^* d \alpha_2 = - d\alpha_1$.
\qed

When the hypersurfaces $\S_1 \subset V$ and $\S_2^* \subset V^*$ are
symmetric about the origin, we can say something more on the above
constructions. For example, in this case we can modify the anti-symplectic
diffeomorphism $\Phi$ in lemma~\ref{symplectic-lemma} to the symplectic
diffeomorphism
\begin{eqnarray*}
\Psi : D^*(\S_1,\S_2^*) & \longrightarrow & D^*(\S_2^*,\S_1) \\
(q,p) & \longmapsto & \Phi(q,-p) \, .
\end{eqnarray*}

Likewise, if we use the diffeomorphism $\phi$ of lemma~\ref{contact-lemma}
to define
\begin{eqnarray*}
\psi : S^*(\S_1,\S_2^*) & \longrightarrow & S^*(\S_2^*,\S_1) \\
(q,p) & \longmapsto & \phi(q,-p) ,
\end{eqnarray*}
then $\psi^* \alpha_1 - \alpha_2$ is an exact $1$-form, and, hence, it
takes closed characteristics of $S^*(\S_1,\S_2^*)$ to closed characteristics
of $S^*(\S_2^*,\S_1)$ while preserving both their orientation and their action.

Moreover, and this will be important in the proof of Sch\"affer's conjecture,
{\it if $\psi(q,p) = (P,Q)$, then $\psi(-q,-p) = (-P,-Q)$.}

\section{Solution of Sch\"affer's conjecture}

Before presenting the proof of Sch\"affer's conjecture, let us derive a
more direct consequence of lemma~\ref{contact-lemma}.

\begin{theorem}\label{contact-equivalence}
Let $\|\cdot\|_{1}$ and $\|\cdot\|_{2}$ be two Minkowski norms on a
finite-dimensional real vector space $V$, and let $\S_{1}$ and $\S_{2}$ denote
their unit spheres. The volume and the length spectrum of the Finsler
metric on $\S_{1}$ induced by its embedding into $(V,\|\cdot\|_{2})$
are equal, respectively, to the volume and the length spectrum of the
Finsler metric on $\S_{2}^{*}$ induced by its embedding into
$(V,\|\cdot\|_{1}^{*})$.
\end{theorem}

\Proof
Notice that the unit co-sphere bundle of $\S_1 \subset (V,\|\cdot\|_2)$ is
$S^*(\S_1,\S_2^*)$ and that the unit co-sphere bundle of
$\S_2^* \subset (V^*,\|\cdot\|_1^*)$ is $S^*(\S_2^*,\S_1)$.
By lemma~\ref{contact-lemma}, these two hypersurfaces are exact
contactomorphic. Therefore, by proposition~\ref{invariants},
their volumes and their action spectra coincide.
\qed

Now we are now ready to prove Sch\"affer's conjecture:

\begin{theorem}
The girth of a normed space equals the girth of its dual.
\end{theorem}

\Proof
We shall need two results of Sch\"affer \cite{Schaffer:2}:
\begin{enumerate}
\item It is enough to prove the conjecture for finite-dimensional normed spaces.
\item The girth is an invariant of isometry classes of normed spaces that is
      continuous with respect to Banach-Mazur topology.
\end{enumerate}

Since (isometry classes of) Minkowski spaces are dense in the Banach-Mazur
compactum of isometry classes of normed spaces with a fixed dimension, then
it is enough for us to prove that {\it the girth of a Minkowski space equals the
girth of its dual.}

If $(V,\|\cdot\|_1)$ and $(V,\|\cdot\|_2)$ are two Minkowski spaces with
unit spheres $\S_1$ and $\S_2$, respectively, we know from the last paragraph
of the previous section that there exists an exact contactomorphism
$\psi : S^*(\S_1,\S_2^*) \rightarrow S^*(\S_2^*,\S_1)$ that takes closed
characteristics to closed characteristics while preserving both their
orientation and action. Recall also that if $\psi(q,p) = (P,Q)$, then
$\psi(-q,-p) = (-P,-Q)$. This means that any characteristic in
$S^*(\S_1,\S_2^*)$ that is invariant under the map $(q,p) \mapsto (-q,-p)$
is sent to a characteristic in $S^*(\S_2^*,\S_1)$ that is invariant under
the map $(P,Q) \mapsto (-P,-Q)$.

Since such characteristics are in one-to-one correspondence to geodesics
that are centrally symmetric, we have that the set of lengths of all
centrally symmetric closed geodesics in $\S_1 \subset (V,\|\cdot\|_2)$ equals the
set of lengths of all centrally symmetric closed geodesics in
$\S_2^* \subset (V^*,\|\cdot\|_1^*)$.

When $\|\cdot\|_1 = \|\cdot\|_2$ this implies that the girth of $(V,\|\cdot\|)$
equals the girth of $(V^*,\|\cdot\|^*)$.
\qed

Notice that lemma~\ref{symplectic-lemma} has the following immediate
consequence:

\begin{theorem}\label{symplectic-equivalence}
Let $(V, \|\cdot\|_{1})$ and $(V, \|\cdot\|_{2})$ be two
Minkowski spaces and let $\S_{1}$ and $\S_{2}$ denote their unit
spheres. The unit co-disc bundle of the Finsler metric on $\S_{1}$
induced by its embedding into $(V,\|\cdot\|_{2})$ is symplectomorphic
to the unit co-disc bundle of the Finsler metric on $\S_{2}^{*}$ induced
by its embedding into $(V,\|\cdot\|_{1}^{*})$.
\end{theorem}

Sch\"affer has given examples of finite-dimensional normed spaces where the
unit sphere and the dual unit sphere have different diameters. These norms can
be approximated by Minkowski norms to yields Minkowski spaces with the same
property. From the theorem above we deduce that {\it the diameter of a Finsler
manifold is not a symplectic invariant of its unit co-disc bundle.}

An interesting billiard version of theorems~\ref{contact-equivalence}
and \ref{symplectic-equivalence} due to S. Tabachnikov and E. Gutkin can be
found in their paper \cite{Tabachnikov-Gutkin}.

\section{Integral geometry}

In this section we uncover the symplectic underpinnings of
El-Ekhtiar's generalization of the Cauchy-Crofton formula to
hypersurfaces in finite-dimensional normed spaces (see
\cite{El-Ekhtiar}).

Our first observation is that if $(V,\|\cdot\|)$ is a Minkowski
space, it is possible to define a natural symplectic structure on
the manifold of oriented lines in $V$, $H_1^+(V)$, which makes it
naturally symplectomorphic to the cotangent bundle of the dual
unit sphere,  $T^* \S^*$.

Consider the diagram
$$
\xymatrix{
      &  V \times \S^* \ar[dl]_{\pi_1}\ar[dr]^{\pi_2} &   \\
     H_1^+(V) &     &  T^* \S^* ,}
$$
where $\pi_1$ is the map that sends a point $(q,P)$ to the oriented line
passing through $q$ in the direction of $\leg(P)$, the Legendre transform
of $P$, and $\pi_2$ is the map $(q,P) \mapsto (P,q_P)$.

If $i : V \times \S^* \rightarrow V \times V^\ast$ is the
natural inclusion, and $\omega$ is the symplectic form on $V \times V^*$,
then define $\omega_1$ to be the unique $2$-form on $H_1^+(V)$ such that
$\pi_1^* \omega_1 = i^* \omega$. This is the natural symplectic form on
the space of oriented lines in $V$ seen as the space of geodesics of the
Finsler manifold $(V,\|\cdot\|)$ (see \cite{Alvarez:thesis}). It is not hard
to see that if $\omega_2$ is the symplectic form on $T^* \S^*$, then
$\pi_2^* \omega_2 = i^* \omega$.

In order to identify $H_1^+(V)$ and $T^* \S^*$ notice that the fibers of the
projections $\pi_1$ and $\pi_2$ coincide. Indeed, $\pi_1(q,P) = \pi_1(q',P')$
if and only if $P = P'$ and $q - q'$ is a multiple of $\leg(P)$ or, equivalently,
if $q - q'$ vanishes on $T_P \S^*$. This is the same condition for
$\pi_2(q,P) = \pi_2(q',P')$. Thus, we obtain a diffeomorphism from $H_1^+(V)$
to $T^* \S^*$ by sending a line $\ell$ to the point $\pi_2(\pi_1^{-1}\{\ell\})$.
The equation $ \pi_1^* \omega_1 = i^* \omega = \pi_{2}^* \omega_2$ implies that
this map is a symplectomorphism.

\begin{theorem}\label{integral-geometry}
Let $(V,\|\cdot\|)$ be a Minkowski  space  and let  $M \subset V$
be a smooth, quadratically convex hypersurface. The unit co-disc
bundle for the induced Finsler metric on $M$ and the set  of all
oriented lines in $V$ which pass through the interior of $M$ are
symplectomorphic.
\end{theorem}

\Proof
Applying lemma~\ref{symplectic-lemma} with $\S_1 := M$ and $\S_2^* = \S^*$,
we have an anti-symplectic diffeomorphism $\Phi$ between the unit co-disc
bundle of the Finsler metric on $M$ induced from its embedding in
$(V,\|\cdot\|)$ and the open set $D^*(\S^*,M) \subset T^* \S^*$.

Composing $\Phi$ with the involution $(q,p) \mapsto (q,-p)$ on
$D^*(M,\S^*) \subset T^*M$ (we are using now that $\S^*$ is symmetric
about the origin), we obtain a symplectomorphism. The theorem now follows
from the fact that the natural symplectomorphism between  $T^*\S^*$ and
$H_{1}^{+}(V)$ takes a point $(P,q_P) \in T^* \S^*$ and sends it to the line
passing through $q$ in the direction of the Legendre transform of $P$.
\qed

If we define the volume density on the manifold of oriented lines of
the $n$-dimensional Minkowski space $V$, $H_1^+(V)$, as $|\omega_1^{n-1}|/(n-1)!$,
then the previous theorem immediately implies the Cauchy-Crofton formula for
finite-dimensional normed spaces.

\begin{corollary}[El-Ekhtiar, \cite{El-Ekhtiar}].
Let $(V,\|\cdot\|)$ be a normed space of dimension $n$ and let  $M
\subset V$ be a convex hypersurface. The volume of the set of
lines passing through $M$ equals the volume of $M$ times the
volume of the Euclidean unit ball of dimension $n-1$.
\end{corollary}

\Proof
If we assume that $M$ and the unit sphere of  $(V,\|\cdot\|)$ are
quadratically convex, then the result follows from the previous
theorem and the definition of volume in the manifold of oriented
lines. An approximation argument takes care of the general case.
\qed

\centerline{\sc Acknowledgements }
\medskip
The author is happy to thank Michael Entov for his comments on an early
version of this paper, and is grateful to Serge Tabachnikov for many
interesting discussions.


\end{document}